\documentclass[12pt,a4paper]{amsart}
\usepackage{amscd}

\newcommand{\g}{\mathfrak{g}}
\newcommand{\fg}{\mathfrak{g}}

\newcommand{\fh}{\mathfrak{h}}

\newcommand{\ft}{\mathfrak{t}}
\newcommand{\fa}{\mathfrak{a}}

\newcommand{\fp}{{\mathfrak{g}_P}}
\newcommand{\fsl}{\mathfrak{sl}}
\newcommand{\fsp}{\mathfrak{sp}}

\newcommand{\fosp}{\mathfrak{osp}}
\newcommand{\fgl}{\mathfrak{gl}}
\newcommand{\fso}{\mathfrak{so}}
\newcommand{\der}{\mathfrak{der}}
\newcommand{\mnt}{\mathfrak{int}}
\newcommand{\tri}{\mathfrak{tri}}

\newcommand{\bZ}{\mathbb{Z}}

\newcommand{\bR}{\mathbb{R}}
\newcommand{\bC}{\mathbb{C}}
\newcommand{\bH}{\mathbb{H}}
\newcommand{\bS}{\mathbb{S}}
\newcommand{\bO}{\mathbb{O}}
\newcommand{\bA}{\mathbb{A}}
\newcommand{\bB}{\mathbb{B}}

\DeclareMathOperator{\Hom}{Hom}
\DeclareMathOperator{\Aut}{Aut}
\DeclareMathOperator{\End}{End}

\DeclareMathOperator{\Img}{Im}
\DeclareMathOperator{\so}{SO}
\DeclareMathOperator{\spl}{SL}

\DeclareMathOperator{\Tri}{Tri}
\DeclareMathOperator{\Int}{Int}

\newtheorem{defn}{Definition}[section]

\begin{document}
\title{Sextonions and the magic square}
\author{Bruce W. Westbury}
\address{Mathematics Institute\\
University of Warwick\\
Coventry CV4 7AL}
\email{bww@maths.warwick.ac.uk}
\date{8 October 2004}
\begin{abstract}
Associated to any complex simple Lie algebra is a non-reductive complex
Lie algebra which we call the intermediate Lie algebra. We propose that
these algebras  can be included in both the magic square and the magic
triangle to give an additional row and column. The extra row and column
in the magic square corresponds to the sextonions. This is a six dimensional
subalgebra of the split octonions which contains the split quaternions.
\end{abstract}
\maketitle

\section{Introduction} The Freudenthal magic square is a $4\times 4$
array of complex semisimple Lie algebras. The rows and columns are
indexed by the real division algebras and the square is symmetric.
This is magic because the row (or column) indexed by the octonions 
consists of four of the five exceptional simple Lie algebras.
There are three constructions which give this square namely the
Tits construction, the Vinberg construction and the triality
construction. Each of these constructions can be extended to give
a rectangle of Lie algebras. There is an alternative point of
view which gives a triangle of Lie algebras.

In this paper we introduce the sextonions as a six dimensional
real alternative algebra intermediate between the split quaternions
and the split octonions. Then we argue that the above magic square,
magic rectangle and magic triangle should all be extended to include
an extra row and column. If the rows or columns are indexed by division
algebras then this extra row or column is indexed by the sextonions.
Here is the extended magic square:

\begin{equation}\label{magic}
\begin{tabular}{|c|c|c|c|c|} \hline
$A_1$ & $A_2$ & $C_3$ & $C_3.H_{14}$ & $F_4$ \\ \hline
$A_2$ & $2A_2$ & $A_5$ & $A_5.H_{20}$ & $E_6$ \\ \hline
$C_3$ & $A_5$ & $D_6$ & $D_6.H_{32}$ & $E_7$ \\ \hline
$C_3.H_{14}$ & $A_5.H_{20}$%
 & $D_6.H_{32}$ & $D_6.H_{32}.H_{44}$ & $E_7.H_{56}$ \\ \hline
$F_4$ & $E_6$ & $E_7$ & $E_7.H_{56}$ & $E_8$ \\ \hline
\end{tabular}
\end{equation}
The notation in this table is that $G.H_n$ means that $G$ has
a representation $V$ of dimension $n$ with an invariant symplectic
form, $\omega$. Then $H_n$ means the Heisenberg algebra of $(V,\omega)$
and $G.H_n$ means the semidirect product of $G$ and $H_n$.

These intermediate algebras are also examples of a more general
construction. For example, the intermediate algebras for the
symplectic algebras are the odd symplectic algebras whose
character theory is studied in \cite{MR89c:20065}, \cite{MR91h:20065},
\cite{MR94k:17011} and \cite{MR2000c:20067}.

Our two general references are the survey article \cite{MR2003f:17003}
on the four real division algebras and \cite{MR2020553} which gives
the three constructions of the magic squares of real Lie algebras
and gives isomorphisms between the Lie algebras given by these 
constructions.

This article is a revised version of the preprint referred to in
\cite{math.RT/0402157}. There is some overlap between these two
articles.

\section{Intermediate Lie algebras}
Our discussion of intermediate algebras is based on the grading
associated to extremal elements. The main application of these
has been to the study of simple modular Lie algebras (see,
for example \cite{MR57:6122}). Another application is in
\cite{MR2001k:17007}.

\begin{defn} A triple in $\fg$ is a set of three elements of $\fg$,
$\{E,F,H\}$ such that
\[ [E,F] = H , [H,E] = 2E , [H,F] = -2F \]
\end{defn}

\begin{defn} An element $e\in\fg$ is extremal if the
one dimensional space with basis $e$ is an inner ideal.
This means that for all $y\in \fg$, $[e,[e,y]]$ is a scalar multiple
of $e$.
A triple $(E,H,F)$ is principal if $E$ (and therefore $F$) is
extremal.
\end{defn}

Let $\fg$ be a complex simple Lie algebra. Then principal
triples can be constructed by choosing a Borel subalgebra and
a root $\alpha$ with the same length as the highest root.
Then there is a principal triple with $E$ in the root space
of $\alpha$ and $F$ in the root space of $-\alpha$ and $H=[E,F]$.

Conversely every principal triple arises this way. Let $\{E,H,F\}$
be a principal triple. Let $\overline{\fg}$ be the centraliser of this triple
and let $\overline{\fh}$ be a Cartan subalgebra of $\overline{\fg}$.
Then a Cartan subalgebra of $\fg$ is given by taking the direct sum of
$\fh$ with the vector space spanned by $H$. Then both $E$ and $F$ span
root spaces; the roots are of the form $\pm\alpha$ and have the same length
as the highest root.

In particular this shows that principal triples are unique up to
automorphism of $\fg$.

Any triple gives a grading on $\fg$ by taking the eigenspaces
of $H$. For a principal triple this grading has the following form
\[ \begin{array}{ccccc}
-2 & -1 & 0 & 1 & 2 \\
\bC & V & \overline{\fg}\oplus \bC & V & \bC
\end{array} \]
\begin{defn}
The graded subalgebra 
\[ \begin{array}{ccc}
0 & 1 & 2 \\
\overline{\fg} & V & \bC
\end{array} \]
is the intermediate subalgebra and is denoted by $\widetilde{\fg}$.
\end{defn}
The vector space $V$ has a $\overline{\fg}$-invariant symplectic form.
The structure of the intermediate algebra is that it is the semidirect 
product of $\overline{\fg}$ and $H_V$ the Heisenberg algebra of $V$.
In particular it is not reductive since the radical is $H_V$.

The intermediate algebra $\widetilde{\fg}$ is also the centraliser of
the extremal vector $e$ in the triple and the Levi subalgebra
$\overline{\fg}$ is the centraliser of the triple.

The $\fsl(2)$ subalgebra with basis the triple is also the centraliser
of $\overline{\fg}$ and these are a dual reductive pair. The restriction
of the adjoint representation of $\fg$ to the subalgebra 
$\overline{\fg}\oplus\fsl(2)$ decomposes as
\begin{equation}\label{dur}
\overline{\fg}\otimes 1\oplus 1\otimes\fsl(2) \oplus V\otimes A 
\end{equation}
where $A$ is the two-dimensional fundamental representation of $\fsl(2)$.
Conversely if we are given an algebra $\fg$ with this structure then
we can take a triple in $\fsl(A)$ as a principal triple.

There are some variations on the definition of the intermediate algebra.
Namely let $\fp$ be
the subalgebra given by taking the subspaces in degrees 0, 1 and 2.
Then this is a parabolic subalgebra. Let $\widetilde{\fg}^\prime$
and $\fp^\prime$ be the respective derived subalgebras. Then we have
a commutative diagram of graded Lie algebras.
\begin{equation}\label{cg}
\begin{CD}
\widetilde{\fg} @>>> \fp \\
@VVV @VVV \\
\widetilde{\fg}^\prime @>>> \fp^\prime \\
\end{CD}
\end{equation}
In this diagram the horizontal arrows are both inclusions of subalgebras
of codimension one. The vertical arrows are both one dimensional
central extensions. In addition we have that $\fp^\prime\cong
\widetilde{\fg}^\prime\oplus \bC$ and so the bottom arrow splits.

The main reason for considering the two algebras $\widetilde{\fg}$ and
$\fp$ is that they are both subalgebras of $\fg$. The main reason for
considering the two algebras $\widetilde{\fg}^\prime$ and $\fp^\prime$ is
that they arise when considering finite dimensional representations.
More precisely, the two vertical arrows in (\ref{cg}) give restriction
functors on the categories of finite dimensional representations.
These two functors are isomorphisms.

\subsection{Examples}
For the special linear algebras this structure can be seen as follows.
Let $U$ and $V$ be any two vector spaces. Then $\fgl(U)\oplus\fgl(V)$
is a subalgebra of $\fg = \fgl(U\oplus V)$. Then the restriction of the
adjoint representation of $\fg$ to this subalgebra decomposes as
\[ \fgl(U)\oplus\fgl(V)\oplus U\otimes V^* \oplus U^*\otimes V \]
If we take the special linear group then we get
\[ \fgl(U)\oplus\fsl(V)\oplus U\otimes V^* \oplus U^*\otimes V \]
If we take $V$ to be two dimensional then $V$ and $V^*$ are equivalent
representations and so we see that
\begin{equation}\label{isl}
\overline{\fsl(n+2)} = \fgl(n)
\end{equation}
and the symplectic representation is the sum of the vector representation
and its dual.

For the symplectic algebras this structure can be seen as follows.
Let $U$ and $V$ be symplectic vector spaces. Then $\fsp(U)\oplus\fsp(V)$
is a subalgebra of $\fg = \fsp(U\oplus V)$. Then the restriction of the
adjoint representation of $\fg$ to this subalgebra decomposes as
\begin{equation}\label{spp}
\fsp(U)\oplus\fsp(V)\oplus U\otimes V
\end{equation}
Taking $V$ to be two dimensional we see that
\[ \overline{\fsp(2n+2)} = \fsp(2n) \]
and the symplectic representation is the vector representation.
These are the Lie algebras of Lie groups known as intermediate symplectic
groups or odd symplectic groups. The characters and representations
of these groups are studied in \cite{MR89c:20065} and \cite{MR94k:17011}.

For the special orthogonal algebras this structure can be seen as follows.
Let $U$ and $V$ be vector spaces with non-degenerate symmetric inner products.
Then $\fso(U)\oplus\fso(V)$
is a subalgebra of $\fg = \fso(U\oplus V)$. Then the restriction of the
adjoint representation of $\fg$ to this subalgebra decomposes as
\begin{equation}\label{sop}
\fso(U)\oplus\fso(V)\oplus U\otimes V
\end{equation}
Taking $V$ to be four dimensional and using the isomorphism 
$\fso(4)\cong\fso(3)\oplus\fso(3)$ we see that
\[ \overline{\fso(n+4)} = \fso(3)\oplus\fso(n) \]
and the symplectic representation is the tensor product of the two
dimensional representation of $\fso(3)$ with the vector representation
of $\fso(n)$.

For the exceptional simple Lie algebras we have the following table
\begin{center}
\begin{tabular}{c|ccccc}
$\fg$ & $G_2$ & $F_4$ & $E_6$ & $E_7$ & $E_8$ \\ \hline
$\overline{\fg}$ & $A_1$ & $C_3$ & $A_5$ & $D_6$ & $E_7$ \\ 
$\dim(V)$ & 4 & 14 & 20 & 32 & 56 \\
\end{tabular}
\end{center}
In all these five cases the representation $V$ is irreducible.

An observation here is that $\overline{\fg}$ can also be described
as the reductive Lie algebra whose rank is one less than the rank
of $\fg$ and where the Dynkin diagram is given by removing the
support of the highest root from the Dynkin diagram of $\fg$.

\subsection{Superalgebras}
This structure can also be extended to the basic Lie superalgebras
by choosing a triple in the even algebra. Our notation for the
dimension of a super space is $(n|m)$ where $n$ is the dimension
of the even subspace and $m$ is the dimension of the odd subspace.

For the special linear algebras this structure can be seen as follows.
Let $U$ and $V$ be any two superspaces. Then $\fgl(U)\oplus\fgl(V)$
is a subalgebra of $\fg = \fgl(U\oplus V)$. Then the restriction of the
adjoint representation of $\fg$ to this subalgebra decomposes as
\[ \fgl(U)\oplus\fgl(V)\oplus U\otimes V^* \oplus U^*\otimes V \]
If we take the special linear group then we get
\begin{equation}\label{supsl}
\fgl(U)\oplus\fsl(V)\oplus U\otimes V^* \oplus U^*\otimes V
\end{equation}
If we take $V$ to have dimension $(2|0)$ then $V$ and $V^*$ are equivalent
representations and so we see that
\[ \overline{\fsl(n+2|m)} = \fgl(n|m) \]
and the symplectic representation is the sum of the vector representation
and its dual. 
Alternatively, if we take $V$ to have dimension $(0|2)$ then $V$ and $V^*$
are equivalent representations and so we see that
\[ \overline{\fsl(n|m+2)} = \fgl(n|m) \]
and the symplectic representation is the sum of the vector representation
and its dual. 

For the orthosympectic algebras this structure can be seen as follows.
Let $U$ and $V$ be superspaces with non-degenerate symmetric inner products.
Then $\fosp(U)\oplus\fosp(V)$
is a subalgebra of $\fg = \fosp(U\oplus V)$. Then the restriction of the
adjoint representation of $\fg$ to this subalgebra decomposes as
\begin{equation}\label{suposp}
\fosp(U)\oplus\fosp(V)\oplus U\otimes V
\end{equation}
This includes (\ref{sop}) and (\ref{spp}) as special cases.
Taking $V$ to have dimension $(4|0)$ and using the isomorphism 
$\fso(4)\cong\fso(3)\oplus\fso(3)$ we see that
\[ \overline{\fosp(n+4|m)} = \fso(3)\oplus\fosp(n|m) \]
and the symplectic representation is the tensor product of the two
dimensional representation of $\fso(3)$ with the vector representation
of $\fso(n)$.
Alternatively, taking $V$ to have dimension $(0|2)$ we see that
\[ \overline{\fosp(n|m+2)} = \fosp(n|m) \]
and the symplectic representation is the vector representation
of $\fosp(n|m)$.

For the exceptional simple Lie superalgebra $G(3)$ the even algebra
is $A_1\oplus G_2$ and the odd part is the tensor product of the
two dimensional fundamental representation of $A_1$ with the seven
dimensional fundamental representation of $G_2$. This means we can
take $\overline{G(3)}=G_2$ and the symplectic representation is the
super space $(0|V)$ where $V$ is the seven dimensional fundamental
representation of $G_2$.

For the exceptional simple Lie superalgebra $F(4)$ the even algebra
is $A_1\oplus B_3$ and the odd part is the tensor product of the
two dimensional fundamental representation of $A_1$ with the eight
dimensional spin representation of $B_3$. This means we can
take $\overline{F(4)}=B_3$ and the symplectic representation is the
super space $(0|V)$ where $V$ is the eight dimensional spin
representation of $B_3$.

\section{Sextonions} In this section we construct the sextonions.
This is a six dimensional real algebra. This is a subalgebra
of the split octonions which is closed under conjugation.
This algebra was explicitly constructed in \cite{MR42:7727}.
This algebra was used in \cite{MR43:2039} to study the conjugacy
classes in $G_2$ in characteristics other than 2 or 3.

The real normed division algebras are the real numbers, the complex
numbers, the quaternions and the octonions. These are denoted by
\[ \bR,\bC,\bH,\bO \]
Each algebra is obtained from the previous one by Cayley-Dickson
doubling. These can be complexified to give complex algebras.
These complex algebras are
\[ \bR\otimes\bC = \bC,
\bC\otimes\bC = \bC\oplus \bC,
\bH\otimes\bC = M_2(\bC),
\bO\otimes\bC \]
The three complex algebras other than $\bC$ have a second real form.
These real forms are denoted $\widetilde{\bC}$,$\widetilde{\bH}$ and
$\widetilde{\bO}$. There are isomorphisms
$\widetilde{\bC} = \bR\oplus\bR$ and $\widetilde{\bH} = M_2(\bR)$.
The normed division
algebras are called the compact forms and this second real form
is called the split real form. These split real forms are composition
algebras but are not division algebras.
The sextonions are intermediate between the split quaternions
and the split octonions.

The sextonions can be constructed as follows. The split quaternions
are isomorphic to the algebra of $2\times 2$ matrices. The norm
is given by the determinant and the conjugate of a matrix is the
adjoint matrix. This algebra has a unique alternative bimodule
which is not associative. This is the two dimensional Cayley
module. This result is given in \cite{MR16:330d}. This bimodule
can be constructed by taking a simple left module $M$ and defining
a right action by
\[ qm=m\Bar q \]
for all $q\in \widetilde{\bH}$ and all $m\in M$.
\begin{defn}\label{sext}
Let $\widetilde{\bS}$ be the split
null extension of $\widetilde{\bH}$ by $M$. This means that we put 
$\widetilde{\bS} = \widetilde{\bH} \oplus M$ and define a multiplication by
\[ (q_1,m_1)(q_2,m_2) = (q_1q_2,q_1m_2+m_1q_2) \]
for all $q_1,q_2\in \widetilde{\bH}$ and all $m_1,m_2\in M$. 
The norm is given by
\[ N(q,m)=\det(q) \] 
and if $x=(q,m)$ then $\Bar x = (\Bar q,-m)$.
\end{defn}
Next we show that this is a subalgebra of the split octonions.
The split octonions can be constructed from the split quaternions
by the Cayley-Dickson doubling process. Put $\widetilde{\bO}
=\widetilde{\bH}\oplus\widetilde{\bH}$ and define a multiplication
by
\begin{equation}\label{oct}
(A_1,B_1)(A_2,B_2)=(A_1A_2-\varepsilon B_2\overline{B_1},
A_1B_2+\overline{A_2}B_1) 
\end{equation}
and define $\overline{(A,B)}=(\overline{A},-B)$ and 
$|(A,B)|=|A|+\varepsilon |B|$.

If we apply this to $\bH$ and take $\varepsilon >0$ then we get
the compact octonions and if $\varepsilon <0$ then we get the
split octonions.
If we apply this to the split quaternions then we get the split
octonions for all $\varepsilon\ne 0$.

Then we see that if we take $B$ to have zero second column then
we obtain the sextonions as a subalgebra.

Note also that we have two commuting actions of $\spl_2(\bR)$.
Let $X\in \spl_2(\bR)$ so $\overline{X}=X^{-1}$. Then these actions
are given by
\begin{equation}
(A,B)\mapsto (XAX^{-1},XB) 
\end{equation}
\begin{equation}\label{acts}
(A,B)\mapsto (A,B\overline{X}) 
\end{equation}

The sextonions are not a division algebra or a composition algebra or a
normed algebra since there is a non-trivial radical given by the Cayley
module $M$ and this is the null space for the inner product.
However they are a subalgebra of the split octonions which
is closed under conjugation. There is a multiplication, 
a conjugation and an inner product which are given in 
Definition (\ref{sext}). This structure is also given by restriction
on the split octonions so any identities which involve this
structure and which hold in the split octonions also hold in the sextonions.

The octonions have a 3-step $\bZ$-grading.
The map $A\mapsto (A,0)$ is an inclusion of
$\widetilde{\bH}$ in $\widetilde{\bO}$ and we take the image to
be the subspace of degree zero. The subspace of pairs of the
form $(0,B)$ is a left module. This has a decomposition as a
left module into a subspace $U_-$ where the second column of $B$
is zero and $U_+$ where the first column of $B$ is zero. Take
$U_-$ to be the subspace of degree -1 and $U_+$ to be the subspace
of degree 1. Note that the product of two elements of $U_-$ or $U_+$
is zero so this is a grading.

\begin{equation}\label{octg}
\begin{tabular}{ccc}
 -1 & 0 & 1 \\ \hline
$U$ & $\widetilde{\bH}$ & $U$ \\
\end{tabular}
\end{equation}

Since the multiplication is needed for later calculations we give
it here explicitly. This is closely related to the description
of the split octonions in \cite{math.RT/0402157}.

\begin{equation}\label{octm}
\left(\begin{array}{c}
u_1 \\ A _1 \\ v_1
\end{array}\right)
\left(\begin{array}{c}
u_2 \\ A _2 \\ v_2
\end{array}\right)
=
\left(\begin{array}{c}
A_1u_2 + \overline{A}_2u_1 \\
A_1A_2 + \overline{(u_2,v_2)}(u_1,v_1) \\
A_1v_2 + \overline{A}_2v_1 
\end{array}\right)
\end{equation}
where $(u,v)$ means put the two column vectors $u$ and $v$
side by side to form a matrix.

\section{Elementary series}
There are three simple constructions which associate a Lie algebra to
the four normed division algebras. These can all be extended to the
sextonions. In this section we show that for each of these constructions
we have that $\overline{\fg(\widetilde{\bO})}=\fg(\widetilde{\bH})$
and that the intermediate algebra is $\fg(\widetilde{\bS})$.

\subsection{Derivations}
The first construction is the derivation algebra. The derivation algebras are
\begin{center}
\begin{tabular}{ccccc}
$\bR$ & $\bC$ & $\bH$ & $\bS$ & $\bO$ \\ \hline
$0$ & $0$ & $A_1$ & $A_1.H_4$ & $G_2$ \\
\end{tabular}
\end{center}

First we look at the derivation algebra of $\widetilde{\bO}$.
The model we take for this is the Cayley-Dickson double of 
$\widetilde{\bH}\cong M_2(\bR)$ given in (\ref{oct}).

The grading in (\ref{octg}) induces a 5-step $\bZ$-grading on the
derivation algebra. This grading is given by
\begin{equation}\label{derg}
\begin{tabular}{ccccc}
-2 & -1 & 0 & 1 & 2 \\ \hline
$\bR$ & $V$ & $\der(\widetilde{\bH})\oplus\bR$ & $V$ & $\bR$ \\
\end{tabular}
\end{equation}
where $V$ is the four dimensional irreducible representation of
$\der(\widetilde{\bH})\cong \fsl_2(\bR)$.

Let $E$ and $F$ be the maps
\[ E \colon (u,A,v) \mapsto (0,0,u) \]
\[ F \colon (u,A,v) \mapsto (v,0,0) \]
Then these can be shown to be derivations by direct calculation.
Put $H=[E,F]$; then $\{E,H,F\}$ is a triple. Note that the 
Lie subalgebra given by this triple is the Lie algebra of the second
action of $\spl_2(\bR)$ in (\ref{acts}).

Also the grading
on $\widetilde{\bO}$ in (\ref{octg}) is also the grading by the
eigenvalues of $H$. This implies that the grading on the derivation
algebra induced by the grading in (\ref{octg}) is also the grading by the
eigenvalues of $H$.

A direct calculation also shows that any derivation of degree two is a
scalar multiple of $E$. Hence the triple $\{E,H,F\}$ is a principal
triple.

Now the derivation algebra of $\bO$ was identified with the Lie
algebra $G_2$ by Elie Cartan in 1915. Take a principal triple
in $G_2$ with $E$ in the highest root space and $F$ in the lowest
root space. Then by inspecting the root diagram we see that the
associated 5-step $\bZ$-grading is given by (\ref{derg}).

Alternatively we can take the construction of the derivation
algebra of $\widetilde{\bO}$ given in \cite{MR96d:22001}.
This construction shows that the derivation algebra has a grading
by the cyclic group of order three with components
$W$, $\spl(W)$, $W^*$ in degrees $-1$, $0$, $1$ where $W$
has dimension three. Then take a principal triple in $\spl(W)$.
The gradings on $W$ and $W^*$ given by the eigenspaces of $H$
are both given by taking one dimensional spaces in dimensions
-1,0 and 1. This shows that the dimensions of the graded components
of the derivation algebra are as given in (\ref{derg}).

Next we consider the derivations of the sextonions. First we consider
$\der^{\widetilde{\bS}}(\widetilde{\bO})$ which consists of the 
derivations of $\widetilde{\bO}$ which preserve $\widetilde{\bS}$.
The restriction homomorphism
$\der^{\widetilde{\bS}}(\widetilde{\bO})\rightarrow\der(\widetilde{\bS})$
is the the homomorphism
$\fp\rightarrow \fp^\prime$ in (\ref{cg})
for $\fg=\der(\widetilde{\bO})$.

There is a homomorphism of graded Lie algebras $\fp^\prime
\rightarrow \der(\widetilde{\bS})$. Our aim now is to show that
this is an isomorphism. It is clear that this is an inclusion
and that both graded Lie algebras have non-zero components only
in degrees zero and one. The graded Lie algebra $\fp^\prime$
has $\fgl(2)$ in degree zero and a four dimensional irreducible
representation in degree one.

The derivations of $\widetilde{\bS}$ of degree zero are a
subspace of $\End(\widetilde{\bH})\oplus\End(U)$. The
derivations in $\End(\widetilde{\bH})$ are the derivations of
$\widetilde{\bH}$ which gives a Lie algebra isomorphic to 
$\fsl(2)$. A calculation shows that a derivation in $\End(U)$
is a scalar multiple of the grading operator $H$.

The derivations of $\widetilde{\bS}$ of degree one are a
subspace of $\Hom(\widetilde{\bH},U)$ which has dimension
eight. This space has an action of the degree zero derivations
and the subspace of derivations is invariant under this action.

A derivation of $\widetilde{\bS}$ of degree one is of the form
\[ (A,u) \mapsto (0,\psi(A)) \]
where $\psi\colon \End(U)\rightarrow U$ is a linear map.
The condition on $\psi$ for this to be a derivation is that
\[ \psi(A_1A_2) = A_1.\psi(A_2)+\overline{A}_2.\psi(A_1) \]
for all $A_1,A_2\in\End(U)$.

Putting $A_2=1$ shows that $\psi(1)=0$.

Now assume that $A_1$ and $A_2$ have zero trace. Then
\[ \psi([A_1,A_2])=2A_1\psi(A_2)-2A_2\psi(A_2) \]
Then if we choose a triple this shows that $\psi(E)$ and
$\psi(F)$ are arbitrary and that these values then determine
$\psi(H)$ and hence $\psi$.

\subsection{Triality}
Let $\bA$ be a composition algebra. Then the triality group
$\Tri(\bA)$ consists of triples $(\theta_1,\theta_2,\theta_3)$
in $\so(\bA)\times\so(\bA)\times\so(\bA)$ such that
\[ \theta_1(a)\theta_2(b)=\theta_2(ab) \]
for all $a,b\in\bA$. Let $\tri(\bA)$ be the Lie algebra
of $\Tri(\bA)$. The triality algebras are
\begin{center}
\begin{tabular}{ccccc}
$\bR$ & $\bC$ & $\bH$ & $\bS$ & $\bO$ \\ \hline
$0$ & $T_2$ & $3A_1$ & $(3A_1).H_8$ & $D_4$ \\
\end{tabular}
\end{center}

The three conditions $\theta_i(1)=1$ define three subgroups.
These three subgroups are isomorphic and any one of them can be
taken as the intermediate group $\Int(\bA)$. Let $\mnt(\bA)$ be the
Lie algebra of $\Int(\bA)$. The intersection of
any two of these intermediate subgroups is the automorphism group,
whose Lie algebra is $\der(\bA)$.

The intermediate algebras are
\begin{center}
\begin{tabular}{ccccc}
$\bR$ & $\bC$ & $\bH$ & $\bS$ & $\bO$ \\ \hline
$0$ & $T_1$ & $2A_1$ & $(2A_1).H_6$ & $B_3$ \\
\end{tabular}
\end{center}

Let $O$ be the orbit of $1\in\bA$ under the action of $\so(\bA)$.
Then we have $\so(\bA)/\Int(\bA)\cong O$ and $\Int(\bA)/\Aut(\bA)\cong O$.
In terms of the Lie algebras we can identify the tangent space of
$1\in O$ with $\Img(\bA)$ and then we have vector space isomorphisms
\[ \fso(\bA) = \mnt(\bA)\oplus \Img(\bA) 
\qquad
 \mnt(\bA) = \der(\bA)\oplus \Img(\bA)
\]

Then the grading on the intermediate algebra $\mnt(\widetilde{\bO})$ is
\begin{equation}\label{intg}
\begin{tabular}{ccccc}
-2 & -1 & 0 & 1 & 2 \\ \hline
$\bR$ & $V\oplus U$ & $\mnt(\widetilde{\bH})\oplus\bR$ & $V\oplus U$ & $\bR$ \\
\end{tabular}
\end{equation}

Then the grading on the triality algebra $\tri(\widetilde{\bO})$ is
\begin{equation}\label{trig}
\begin{tabular}{ccccc}
-2 & -1 & 0 & 1 & 2 \\ \hline
$\bR$ & $V\oplus 2U$ & $\tri(\widetilde{\bH})\oplus\bR$ & $V\oplus 2U$ & $\bR$ \\
\end{tabular}
\end{equation}

\subsection{Superalgebras}
There are also two constructions of Lie superalgebras.
These constructions are given in \cite{MR85c:17011} and \cite{MR2031706}
and \cite{MR2004a:17003}. One construction
is to take $\fg(\bA)$ to be the superspace with even part
$\fsl(A)\oplus \der(\bA)$ and odd part $A\otimes\Img(\bA)$
where $A$ is a two dimensional vector space.
This construction gives the Lie superalgebras
\begin{center}
\begin{tabular}{ccccc}
$\bR$ & $\bC$ & $\bH$ & $\bS$ & $\bO$ \\ \hline
$A_1$ & $B(0,1)$ & $B(1,1)$ & $B(1,1).H_{(4|4)}$ & $G(3)$ \\
\end{tabular}
\end{center}
The grading on $\fg(\widetilde{\bO})$ is
\begin{center}
\begin{tabular}{ccccc}
-2 & -1 & 0 & 1 & 2 \\ \hline
$\bR$ & $(V|U\otimes A)$ & $\fg(\widetilde{\bH})\oplus\bR$%
& $(V|U\otimes A)$ & $\bR$ \\
\end{tabular}
\end{center}

A second construction 
is to take $\fg(\bA)$ to be the superspace with even part
$\fsl(A)\oplus \mnt(\bA)$ and odd part $A\otimes\bA$.
where $A$ is a two dimensional vector space.
This construction gives the Lie superalgebras
\begin{center}
\begin{tabular}{ccccc}
$\bR$ & $\bC$ & $\bH$ & $\bS$ & $\bO$ \\ \hline
$B(0,1)$ & $A(1,0)$ & $D(2,1;\mu)$ & $D(2,1).H_{(6|4)}$ & $F(4)$ \\
\end{tabular}
\end{center}
The grading on $\fg(\widetilde{\bO})$ is
\begin{center}
\begin{tabular}{ccccc}
-2 & -1 & 0 & 1 & 2 \\ \hline
$\bR$ & $(V\oplus U|U\otimes A)$ & $\fg(\widetilde{\bH})\oplus\bR$%
& $(V\oplus U|U\otimes A)$ & $\bR$ \\
\end{tabular}
\end{center}
where we have used the grading (\ref{intg}).

\section{The magic square}
There are three constructions of the magic square. All three constructions
take a pair of composition algebras $(\bA,\bB)$ and produce a semisimple
Lie algebra $L(\bA,\bB)$. The original construction is due to Freudenthal-Tits.
Other constructions are
Vinberg construction and the triality construction. In all these cases we can
extend the construction to include the sextonions and all constructions give
isomorphic Lie algebras. Again we find that the intermediate subalgebra of
$L(\bA,\widetilde{\bO})$ is $L(\bA,\widetilde{\bS})$ and the Levi subalgebra
is $L(\bA,\widetilde{\bH})$.

Let $\bA$ be a composition algebra and $J$ a Jordan algebra.
The Tits construction is
\[ T(\bA,J) = \der(\bA)\oplus \der(J) \oplus \Img(\bA)\otimes\Img(J) \]
Then the grading on $T(\widetilde{\bO},J)$ is
\begin{center}
\begin{tabular}{ccccc}
-2 & -1 & 0 & 1 & 2 \\ \hline
$\bR$ & $V\oplus U\otimes\Img(J)$ & $T(\widetilde{\bH},J)\oplus\bR$ %
& $V\oplus U\otimes\Img(J)$ & $\bR$ \\
\end{tabular}
\end{center}
For the construction of the magic square we take the Jordan algebra
$J=J(\bB)$ to be $H_3(\bB)$ which consists of $3\times 3$ Hermitian
matrices with entries in $\bB$. We can also take $J=0$ which gives the
derivation algebras.

The Vinberg construction is
\[ V(\bA,\bB) = \der(\bA)\oplus\der(\bB) \oplus A_3^\prime(\bA\otimes\bB) \]
where $A_3^\prime(\bA)$ means trace-free anti-Hermitian $3\times 3$
matrices with entries in $\bA$.
Then the grading on $V(\widetilde{\bO},\bB)$ is
\begin{center}
\begin{tabular}{ccccc}
-2 & -1 & 0 & 1 & 2 \\ \hline
$\bR$ & $V\oplus A_3^\prime(U\otimes\bB)$ & $V(\widetilde{\bH},\bB)\oplus\bR$ %
& $V\oplus A_3^\prime(U\otimes\bB)$ & $\bR$ \\
\end{tabular}
\end{center}
Since $U$ is imaginary we can identify $A_3^\prime(U\otimes\bB)$ with
$U\otimes H_3^\prime(\bB)$ where $H_3^\prime(\bA)$ means trace-free
Hermitian $3\times 3$ matrices with entries in $\bA$. This is also
$\Img(J)$ for $J=H_3(\bA)$.

The triality construction is
\[ A(\bA,\bB) = \tri(\bA)\oplus\tri(\bB) \oplus 3(\bA\otimes \bB) \]
Then the grading on $A(\widetilde{\bO},\bB)$ is
\begin{center}
\begin{tabular}{ccccc}
-2 & -1 & 0 & 1 & 2 \\ \hline
$\bR$ & $V\oplus 2U \oplus 3U\otimes\bB$ & $A(\widetilde{\bH},\bB)\oplus\bR$ %
& $V\oplus 2U \oplus 3U\otimes\bB$ & $\bR$ \\
\end{tabular}
\end{center}
where we have used the grading (\ref{trig}).

Taking both of the algebras in the triality construction to be the
sextonions gives a Lie algebra of dimension 144. Let $V_{56}$ be the 56
dimensional fundamental representation of $E_7$. Then the grading
on $E_8$ has components
\[ \begin{array}{ccccc}
-2 & -1 & 0 & 1 & 2 \\
\bC & V_{56} & E_7\oplus \bC & V_{56} & \bC
\end{array} \]
Then take a principal triple in $E_7$.
Then this triple commutes with the principal triple in $E_8$
and so we have a bigrading on $E_8$ with components
\begin{equation}\label{bigrad}
\begin{array}{ccccc}
 & & \bC & & \\
 & V_{12} & S_{32} & V_{12} & \\
\bC & S_{32} & D_6\oplus\bC\oplus\bC & S_{32} & \bC \\
 & V_{12} & S_{32} & V_{12} & \\
 & & \bC & & 
\end{array}
\end{equation}
where $S_{32}$ is a spin representation of $D_6$ of dimension 32
and $V_{12}$ is the vector representation of dimension 12. This constructs
the Lie algebra $E_8$ as
\[ (D_6\oplus\fsl(A)\oplus\fsl(B)) 
\oplus (V\otimes A\otimes B) \oplus S\otimes A \oplus S\otimes B \]
where $A$ and $B$ are two dimensional vector spaces.

Also if we take the total grading in (\ref{bigrad}) we get
the grading with components
\[ \begin{array}{ccccc}
-2 & -1 & 0 & 1 & 2 \\
V_{14} & S_{64} & D_7\oplus \bC & S_{64}  & V_{14}
\end{array} \]
where $S_{64}$ is a spin representation of $D_7$ of dimension 64
and $V_{14}$ is the vector representation of dimension 14. The non-negative
part of this grading gives a second maximal parabolic subgroup of $E_8$.
The even part of this grading is isomorphic to $D_8$. The odd part
is a spin representation of $D_8$. This is used in \cite{MR98b:22001}
to construct the Lie algebra $E_8$.

The bigrading in (\ref{bigrad}) gives the following non-negatively bigraded
Lie algebra of dimension 144.
\[ \begin{array}{ccc}
D_6 & S_{32} & \bC \\
S_{32} & V_{12} & \\
\bC & &
\end{array} \]

\subsection{Exceptional series}
In this section we consider the exceptional series introduced
in \cite{MR96m:22012} and the subexceptional series. These are
finite series of semisimple Lie algebras. The exceptional series
includes all five exceptional simple Lie algebras.
Here we take these Lie algebras to be parametrised by $m$.
Different authors have used other parameters such as the dual
Coxeter number. All of these parameters are related to $m$
by Mobius transformations. 

For the Lie algebras in the magic square we get $L(\bO,\bA)$
in the exceptional series with $m=\dim(\bA)$ and $L(\bH,\bA)$
in the subexceptional series again with $m=\dim(\bA)$. This
gives the last three rows of (\ref{magic}) with columns
labelled by $m=1,2,4,6,8$. The exceptional series also includes
further columns. Four of these columns are given below:

\begin{center}
\begin{tabular}{|c||c|c|c|c|} \hline
$m$ & $-4/3$ & $-1$ & $-2/3$ & $0$ \\
\hline\hline
$\overline{\fg}$ & $0$ & $T$ & $A_1$ & $3A_1$ \\ \hline
$\widetilde{\fg}$ & $0$ & $T.H_2$ & $A_1.H_4$ & $(3A_1).H_8$ \\ \hline
$\fg$ & $A_1$ & $A_2$ & $G_2$ & $D_4$ \\ \hline
\end{tabular}
\end{center}

The column with $m=0$ contains the triality algebras and the
column with $m=-2/3$ contains the derivation algebras.

In this section we extend the exceptional series to include
some simple Lie superalgebras. Let $\fg(\bH)$ be a Lie algebra
in the subexceptional series and $\fg(\bO)$ the corresponding 
Lie algebra in the exceptional series. Then $\fg(\bH)$ has a
distinguished representation $V$ of dimension $6m+8$ which has
a $\fg(\bH)$-invariant symplectic form. This is the representation
$V$ in (\ref{dur}).

This is consistent with the dimension formulae:
\[ \dim(\fg(\bH))=3\frac{(2m+3)(3m+4)}{(m+4)},
 \dim(\fg(\bO))=2\frac{(3m+7)(5m+8)}{(m+4)}
 \]
In these notes we show that this construction also makes sense
for some values of $m$ for which $6m+8$ is a negative integer.
In this case we take $V$ to be an odd superspace and apply the
same construction to obtain a Lie superalgebra.

\begin{center}
\begin{tabular}{|c|cccccc|}
\hline
$m$  & -3 & -8/3 & -5/2 & -7/3 &  -2 & -3/2  \\
\hline
$6m+8$ & -10 & -8 & -7 & -6 & -4 & -1 \\
\hline
$\fg(\bH)$ & $D_5$ & $B_3$ & $G_2$ & $A_2+T$ & $A_1$ & 0 \\
& $\fso(10)$ & $\fso(7)$ & $G_2$ & $\fgl(3)$ & $\fsl(2)$ & 0 \\
\hline
$\fg(\bO)$ & $D(5,1)$ & $F(4)$ & $G(3)$ & $A(2,1)$ & $A(1,1)$ &
$B(1,1)$ \\
 & $\fosp(10|2)$ & $F(4)$ & $G(3)$ & $\fsl(3|2)$ & $\fsl(2|2)$ &
 $\fosp(1|2)$ \\
\hline
\end{tabular}
\end{center}

There is a distinguished representation $V$ of $\fg(\bH)$ dimension $-6m-8$.
The structure that these representations have in common is that
\[ S^2(V) = 1\oplus V^2 \qquad \Lambda^2(V) = \fg\oplus V_2 \]

The representation $V_2$ is somewhat degenerate:
\begin{enumerate}
\item For $m=-3$, $\fg(\bH)=\fso(10)$, $V$ is the vector representation
and $V_2=0$.
\item For $m=-8/3$, $\fg(\bH)=\fso(7)$, $V$ is the spin representation
and $V_2$ is the vector representation.
\item For $m=-5/2$, $\fg(\bH)=G_2$ and $V$ and $V_2$ are both the seven
dimensional fundamental representation.
\item For $m=-7/3$, $\fg(\bH)=\fgl(3)$, $V$ is the sum of the vector
representation and its dual and $V_2$ is the adjoint representation.
\item For $m=-3/2$, $\fg(\bH)=0$, $V$ has dimension one and $V_2=0$.
\end{enumerate}

Note that in some cases we can replace $\fg(\bH)$ by a Lie superalgebra
and still keeping this structure.

\begin{center}
\begin{tabular}{|c|ccc|}
\hline
$m$  & -3 &  -7/3 &  -2 \\
\hline
$6m+8$ & -10 &  -6 & -4 \\
\hline
$\fg(\bH)$ & $D(n+5,n)$ & $A(n+2,n)+T$ & $A(n+1,n)$ \\
& $\fosp(2n+10|2n)$ & $\fgl(n+3|n)$ & $\fsl(n+2|n)$ \\
\hline
$\fg(\bO)$ & $D(n+5,n+1)$ & $A(n+2,n+1)$ & $A(n+1,n+1)$ \\
 & $\fosp(2n+10|2n+2)$ & $\fsl(n+3|n+2)$ & $\fsl(n+2|n+2)$ \\
\hline
\end{tabular}
\end{center}

\begin{center}
\begin{tabular}{|c|ccc|}
\hline
$m$  & -3/2 & -4/3 & -1 \\
\hline
$6m+8$ & -1 & 0 & 2 \\
\hline
$\fg(\bH)$ & $B(n+1,n)$ & $A(n,n)+T$ & $A(n+1,1)+T$ \\
& $\fosp(2n+1|2n)$ & $\fgl(n|n)$ & $\fgl(n+1|1)$ \\
\hline
$\fg(\bO)$ & $B(n+1,n+1)$ & $A(n+1,n)$ & $A(n+2,2)$ \\
 & $\fosp(2n+1|2n+2)$ & $\fsl(n+2|n)$ & $\fsl(n+3|n)$ \\
\hline
\end{tabular}
\end{center}
These follow from the general decompositions in (\ref{supsl})
and (\ref{suposp}).

The point $m=-8/5$ on the exceptional line corresponds to the trivial
Lie algebra. However there is no corresponding Lie algebra on the
subexceptional line.

\subsection{Magic triangle}
There is another approach to the magic square based on dual reductive
pairs. This constructs a magic triangle. This magic triangle is given
in \cite{cvit3}, \cite{MR98j:22013} and \cite{deligne}. This is also
implicit in \cite{julia}.

The involution which sends $\g$ to the centraliser in $E_8$ corresponds
to the involution
\[ m \mapsto \frac{-2m}{m+2} \]
If we include the Lie algebra $E_7.H_{56}$ with $m=6$ then this 
suggests that we should also include a Lie algebra for $m=-3/2$.
This Lie algebra is given as the Lie superalgebra $\fosp(1|2)$.
Taken literally this suggests that $\fosp(1|2)$ and $E_7.H_{56}$
are a dual reductive pair in $E_8$. However $\fosp(1|2)$ is not a
subalgebra and $E_7.H_{56}$ is not reductive.

More generally the decomposition (\ref{dur}) shows that $A_1$ and
$\overline{\fg}$ are a dual reductive pair in $\fg$. Here we do
a formal calculation which shows that as characters of
$\overline{\fg}\oplus \fsl(2)$ we have
\begin{equation}\label{dus}
\fg  = \widetilde{\fg}\otimes 1\oplus 1\otimes\fosp(1|2) \oplus
(V\oplus 1)\otimes A 
\end{equation}
where $A$ is the vector representation of dimension $(2|1)$.

Then we write a super vector space as $V_+-V_-$ where $V_+$ is the even
part and $V_-$ is the odd part. We write $[n]$ for the irreducible highest 
weight representation of $\fsp(1)$ with highest weight $n$ (and dimension $n+1$)
and we regard a representation of $\fosp(1|2)$ as a super representation of
$\fsp(1)$.
In particular, the adjoint representation of $\fosp(1|2)$ is written as
$[2]-[1]$ and the representation $A$ is written as $[1]-[0]$.
Then the right hand side of (\ref{dus}) is
\[
(\overline{\fg}+V+1)\otimes [0] \oplus 1\otimes([2]-[1])
\oplus (V+1)\otimes ([1]-[0])
\]
Expanding this and cancelling equal terms with oposite signs leaves
\[ \overline{\fg}\otimes [0] + 1\otimes [2] + V\otimes [1] \]
which is (\ref{dur}).

If we apply this to the Lie algebras in the exceptional series then
this calculation is our justification for including an extra row
and column in the magic triangle.

\section{Adams series}
The triality construction constructs the Lie algebra $L(\bA,\bB)$
with a $\bZ_2\times \bZ_2$ grading. If we take any one of the three
$\bZ_2$ gradings then in degree zero we get the Lie algebra
\[ \ft(\bA,\bB) = \tri(\bA)\oplus \tri(\bB) \oplus \bA\otimes\bB \]
and in degree one we get the spin representation
$(\bA\otimes \bB) \oplus (\bA\otimes \bB)$.
The table for these Lie algebras is given in \cite{triality}.
Note that $\ft(\bA,\bB)$ is a subalgebra of equal rank in 
$L(\bA,\bB)$.

Comparing this construction and (\ref{sop}) we observe that there is
a variation on these two constructions.
Let $V$ be a vector space with a non-degenerate symmetric inner product.
Then the algebra $\fa(\bA,V)$ is defined as a vector space by
\[ \fa(\bA,V)=\tri(\bA)\oplus\fso(V) \oplus \bA\otimes V \]
The Lie bracket is defined so that $\ft(\bA)\oplus\fso(V)$ is a subalgebra
and $\bA\otimes V$ is the obvious representation. The Lie bracket of two
elements of $\bA\otimes V$ is the usual Lie bracket so that
$\fso(\bA\oplus V)$ is a subalgebra of $\fa(\bA,V)$.

Note that we have inclusions
$\fso(\bA\oplus W)\subset\fa(\bA,W)$ and
$\fa(\bA ,\bB)\subset\ft(\bA,\bB)$.

The grading on $\fa(\widetilde{\bO},W)$ is
\begin{center}
\begin{tabular}{ccccc}
-2 & -1 & 0 & 1 & 2 \\ \hline
$\bR$ & $V\oplus 2U \oplus U\otimes W$ & $\fa(\widetilde{\bH},W)\oplus\bR$ %
& $V\oplus 2U \oplus U\otimes W$ & $\bR$ \\
\end{tabular}
\end{center}
where we have used the grading (\ref{trig}).

The grading on $\ft(\widetilde{\bO},\bA)$ is
\begin{center}
\begin{tabular}{ccccc}
-2 & -1 & 0 & 1 & 2 \\ \hline
$\bR$ & $V\oplus 2U \oplus U\otimes \bA$ & $\ft(\widetilde{\bH},\bA)\oplus\bR$ %
& $V\oplus 2U \oplus U\otimes \bA$ & $\bR$ \\
\end{tabular}
\end{center}
where we have used the grading (\ref{trig}).


This gives the following generalisation of (\ref{bigrad}).
Take a principal triple in $L(\bA,\bO)$ with centraliser
$L(\bA,\bH)$ and then take a principal triple in $L(\bA,\bH)$.
These two triples commute and so we get a bigrading on 
$L(\bA,\bO)$. Put $m=\dim(\bA)$ then this bigrading is given by
\begin{equation}\label{bigrada}
\begin{array}{ccccc}
 & & \bC & & \\
 & V_{m+4} & S_{4m} & V_{m+4} & \\
\bC & S_{4m} & \fa(\bA,W_4)\oplus\bC\oplus\bC & S_{4m} & \bC \\
 & V_{m+4} & S_{4m} & V_{m+4} & \\
 & & \bC & & 
\end{array}
\end{equation}
where $S_{4m}$ is a spin representation of $\fa(\bA,W_4)$ of dimension $4m$
and $V_{m+4}$ is the vector representation of dimension $m+4$. This constructs
the Lie algebra $L(\bA,\bO)$ as
\[ (\fa(\bA,W_4)\oplus\fsl(A)\oplus\fsl(B)) 
\oplus (V\otimes A\otimes B) \oplus S\otimes A \oplus S\otimes B \]
where $A$ and $B$ are two dimensional vector spaces.

Also if we take the total grading in (\ref{bigrada}) we get
the grading with components
\[ \begin{array}{ccccc}
-2 & -1 & 0 & 1 & 2 \\
V_{m+6} & S_{8m} & \fa(\bA,W_6)\oplus \bC & S_{8m}  & V_{m+6}
\end{array} \]
where $S_{8m}$ is a spin representation of $\fa(\bA,W_6)$ of dimension $8m$
and $V_{m+6}$ is the vector representation of dimension $m+6$. The non-negative
part of this grading gives a second maximal parabolic subgroup of $L(\bA,\bO)$.
The even part of this grading is isomorphic to $\fa(\bA,W_8)\cong\ft(\bA,\bO)$.
The odd part is a spin representation of dimension $16m$.
This is used in \cite{MR98b:22001}
to construct the Lie algebra $L(\bA,\bO)$.

\bibliographystyle{halpha}
\bibliography{sex}

\end{document}